\documentclass[10pt,twoside]{siamart1116}

\usepackage[english]{babel}
\usepackage{graphicx,epstopdf,epsfig}
\usepackage{amsfonts,epsfig,fancyhdr,graphics, hyperref,amsmath,amssymb,biblatex}
\newcommand{\Names}{Muhammad Firmansyah Kasim}
\newcommand{\Title}{Derivatives of partial eigendecomposition of a real symmetric matrix for degenerate cases}

\renewcommand{\theequation}{\thesection.\arabic{equation}}
\renewtheorem{theorem}{Theorem}[section]

\begin{document}

\setcounter{page}{1}

\thispagestyle{empty}

 \title{\Title
 \thanks{Date: \today}
 }

\author{
Muhammad F.\ Kasim\thanks{Department of Physics,
University of Oxford, Parks Rd, Oxford OX1 3PU, United Kingdom
(muhammad.kasim@physics.ox.ac.uk). Supported by a generous grant from EPSRC: EP/P015794/1.}
}

\markboth{\Names}{\Title}

\maketitle

\begin{abstract}
This paper presents the forward and backward derivatives of partial eigendecomposition, i.e. where it only obtains some of the eigenpairs, of a real symmetric matrix for degenerate cases.
The numerical calculation of forward and backward derivatives can be implemented even if the degeneracy never disappears and only some eigenpairs are available.
\end{abstract}

\begin{keywords}
Eigendecomposition, forward derivative, backward derivative, degenerate
\end{keywords}
\begin{AMS}
15A18, 65F15
\end{AMS}

\section{Introduction} \label{sec:intro}

In many physics and engineering simulations, eigendecomposition of a real symmetric matrix is a critical step in computing physical quantities of the simulations.
For some applications, e.g. plane-wave density functional theory \cite{gonze2009abinit}, the matrix is too large to be stored in a computer's memory, therefore it can only be represented by its matrix-vector multiplication.
The eigendecomposition of this kind of matrix is usually performed by only taking the $k$-largest or smallest eigenvalues and the corresponding eigenvectors, instead of computing all the eigenpairs.

Although algorithms to compute the partial eigendecomposition are available \cite{knyazev2001lobpcg}, the derivatives (forward and backward) of the partial eigendecomposition is not available, hindering its use in automatic differentiation programs \cite{abadi2016tensorflow, paszke2019pytorch}.
Giles \cite{giles2008extended} presented the forward and backward derivatives of complete eigendecomposition for explicit matrix, but nothing on partial eigendecomposition and it does not explain handling cases with repeated eigenvalues (or known as degenerate cases).
The algorithm to compute the forward derivatives of complete eigendecomposition for degenerate cases was presented by van der Aa, \textit{et al.} \cite{van2007computation} by iteratively computing higher order derivatives until the degeneracy disappear.
However, the work does not explain the backward derivative and cases where the degeneracy never disappears.
It also needs the complete eigenvalues and eigenvectors to implement.

In this work, the forward and backward derivatives of partial eigendecomposition is derived and presented.
The expressions can be implemented numerically even if the degeneracy never disappears and only $k$ eigenpairs are available.

\section{Problem statement and results summary}\label{sec:problem-results}

Let $\mathbf{A}, \mathbf{M}\in\mathbb{R}^{n\times n}$ be real symmetric matrices with $\mathbf{M}$ being positive definite.
The partial eigendecomposition can be written as
\begin{equation}
    \label{eq:orig-eidecomp}
    \mathbf{AX} = \mathbf{MX\Lambda}
\end{equation}
where $\mathbf{X}\in\mathbb{R}^{n\times k}$ is a matrix containing $k$ eigenvectors in its columns ($k < n$), and $\mathbf{\Lambda}\in\mathbb{R}^{k\times k}$ is a diagonal matrix of the $k$ eigenvalues.
The eigenvectors are normalized so that
\begin{equation}
    \mathbf{X}^T\mathbf{MX} = \mathbf{I}.
\end{equation}

\subsection{Results on forward derivative}

If the matrices $\mathbf{A}$ and $\mathbf{M}$ are perturbed by $\mathbf{A'}$ and $\mathbf{M'}$, respectively, then the perturbation of the eigenvalues and eigenvectors are respectively given by
\begin{align}
    \mathbf{\Lambda}' &= \mathbf{I}\circ \left[\mathbf{X}^T \left(\mathbf{A}' \mathbf{X} - \mathbf{M'X\Lambda}\right)\right] \\
    \label{eq:summary-fwdderiv-eivecs}
    \mathbf{X}' &= -\frac{1}{2}\mathbf{X}\left[\mathbf{I}\circ\left(\mathbf{X}^T\mathbf{M'X}\right)\right] - \mathbf{Y}' + \mathbf{X}\left[\mathbf{D}\circ\left(\mathbf{X}^T\mathbf{MY'}\right)\right]
\end{align}
with $\circ$ denotes element-wise multiplication, $\mathbf{I}$ as an $n\times n$ identity matrix, $\mathbf{Y'}$ and $\mathbf{V'}$ are given by
\begin{align}
    \label{eq:summary-fwdderiv-temp-y}
    \mathbf{AY' - MY'\Lambda} &= \mathbf{V}' - \mathbf{MX}\left[\mathbf{I}\circ\left(\mathbf{X}^T\mathbf{V'}\right)\right] \\
    \label{eq:temp-v-fwdderiv}
    \mathbf{V'} &= \mathbf{A'X} - \mathbf{M'X\Lambda},
\end{align}
and $\mathbf{D}$ as the degeneracy matrix with elements
\begin{equation}
    \label{eq:degen-matrix}
    D_{ij} = \begin{cases}
        1, &\ \text{if }\Lambda_{ii} = \Lambda_{jj}\\
        0, &\ \text{otherwise}.
    \end{cases}
\end{equation}
The equation \ref{eq:summary-fwdderiv-temp-y} is Sylvester equation and can be solved numerically.
If matrix $\mathbf{A}$ and $\mathbf{M}$ are only represented implicitly, i.e. only the matrix-vector product is known, then it can be solved using batched iterative linear equation solvers, such as GMRES \cite{saad1986gmres} or MINRES \cite{paige1975minres}.

The eigenvectors perturbation in equation \ref{eq:summary-fwdderiv-eivecs} is only valid if the following condition is satisfied,
\begin{equation}
    \left(\mathbf{D} - \mathbf{I}\right)\circ\left[\mathbf{X}^T(\mathbf{A'X} - \mathbf{M'X\Lambda})\right] = \mathbf{0}.
\end{equation}

\subsection{Results on backward derivative}
Backward derivative would make sense if used in a context when the eigenvalues $\mathbf{\Lambda}$ and eigenvectors $\mathbf{X}$ are used for calculating a value, $\mathcal{L}$.
Let's denote the sensitivity of the value $\mathcal{L}$ with respect to the eigenvalues and eigenvectors as $\mathbf{\overline{\Lambda}} \equiv \partial\mathcal{L}/\partial\mathbf{\Lambda}$ and $\mathbf{\overline{X}} \equiv \partial\mathcal{L}/\partial\mathbf{X}$, respectively.
The sensitivity with respect to the matrices $\mathbf{A}$ and $\mathbf{M}$ respectively are
\begin{align}
    \mathbf{\overline{A}} &= \mathbf{X\overline{\Lambda}X}^T - \mathbf{\overline{V}X}^T\\
    \mathbf{\overline{M}} &= \mathbf{X\Lambda\overline{\Lambda}X}^T -
        \frac{1}{2}\mathbf{X}\left[\mathbf{I}\circ\left(\mathbf{X}^T\mathbf{\overline{X}}\right)\right]\mathbf{X}^T +
        \mathbf{\overline{V}\Lambda X}^T.
\end{align}
where
\begin{align}
    \mathbf{\overline{V}} &= \mathbf{\overline{Y}} -
        \mathbf{X}\left[\mathbf{I}\circ\left(\mathbf{X}^T\mathbf{M\overline{Y}}\right)\right] \\
    \mathbf{A\overline{Y}} - \mathbf{M\overline{Y}E} &=
        \mathbf{\overline{X}} -
        \mathbf{MX} \left[\mathbf{D}\circ\left(\mathbf{X}^T \mathbf{\overline{X}}\right)\right].
\end{align}

If the matrix $\mathbf{A}$ and $\mathbf{M}$ are always symmetric, the expressions above are only valid if
\begin{align}
    (\mathbf{D} - \mathbf{I}) \circ \left(\mathbf{X}^T\mathbf{\overline{X}} - \mathbf{\overline{X}}^T\mathbf{X}\right) &= \mathbf{0} \\
    \nonumber&\mathrm{or}\\
    (\mathbf{D} - \mathbf{I})\circ\left(\mathbf{X}^T\mathbf{A'X}\right) &= \mathbf{0} \\
    (\mathbf{D} - \mathbf{I})\circ\left(\mathbf{X}^T\mathbf{M'X\Lambda}\right) &= \mathbf{0}.
\end{align}

The next sections will present the derivation of the expressions above, starting from the non-degenerate case, then move to the degenerate case.

\section{Derivation for the non-degenerate case}
\subsection{Forward derivative}

To simplify the work, consider only the $j$-th eigenvalue and eigenvector of matrix $\mathbf{A}$,
\begin{equation}\label{eq:single-eigdecomp}
    \mathbf{Ax}_j = \lambda_j\mathbf{Mx}_j,
\end{equation}
where the eigenvector is normalized
\begin{equation}\label{eq:single-norm}
    \mathbf{x}_j^T \mathbf{Mx}_j = 1.
\end{equation}
For simplicity, from now the index $j$ is not written.

Perturbing the equations \ref{eq:single-eigdecomp} and \ref{eq:single-norm} once yields
\begin{align}
    \label{eq:diff-single-eigdecomp}
    \mathbf{A'x} + \mathbf{Ax}' &= \lambda'\mathbf{Mx} + \lambda \mathbf{M'x} + \lambda \mathbf{Mx}' \\
    \label{eq:diff-single-norm}
    \mathbf{x}^T\mathbf{M'x} + 2 \mathbf{x}^T\mathbf{Mx}' &= 0.
\end{align}
The primed variables indicate a small perturbation of the corresponding variable.
Applying $\mathbf{x}^T$ from left on both sides of equation \ref{eq:diff-single-eigdecomp} gives
\begin{equation}
    \mathbf{x}^T\mathbf{A'x} + \mathbf{x}^T\mathbf{Ax}' = \lambda'\mathbf{x}^T\mathbf{Mx} + \lambda \mathbf{x}^T\mathbf{M'x} + \lambda \mathbf{x}^T\mathbf{Mx}'
\end{equation}
By substituting $\mathbf{x}^T\mathbf{Mx} = 1$ from equation \ref{eq:single-norm} and $\mathbf{x}^T\mathbf{A}$ from equation \ref{eq:single-eigdecomp} to the equation above, the forward derivative of the eigenvalue is obtained as
\begin{equation}\label{eq:fwdderiv-single-eival}
    \lambda' = \mathbf{x}^T(\mathbf{A'} - \lambda\mathbf{M'})\mathbf{x}.
\end{equation}

To obtain the perturbation of the eigenvector $\mathbf{x}$, substitute \ref{eq:fwdderiv-single-eival} to \ref{eq:diff-single-eigdecomp} and rearrange it to produce
\begin{equation}
    \label{eq:fwdderiv-before-splitting}
    (\mathbf{A} - \lambda\mathbf{M})\mathbf{x}' = -(\mathbf{I} - \mathbf{Mxx}^T)(\mathbf{A}' - \lambda\mathbf{M}')\mathbf{x}.
\end{equation}
For non-degenerate case, the matrix $(\mathbf{A} - \lambda\mathbf{M})$ has the rank of $n-1$ and nullify all the components parallel to $\mathbf{x}$.
Therefore, the parallel and orthogonal components of $\mathbf{x}'$ with respect to $\mathbf{x}$ need to be treated separately.
Let's denote the parallel and orthogonal components of $\mathbf{x}'$ as $\mathbf{x}_{\parallel}'$ and $\mathbf{x}_{\perp}'$ respectively.
The term $\mathbf{x}'$ can then be written as
\begin{equation}
    \mathbf{x}' = \mathbf{x}_{\parallel}' + \mathbf{x}_{\perp}'
\end{equation}
with properties
\begin{align}
    \label{eq:ortho-property}
    \mathbf{x}_{\perp}' &= (\mathbf{I} - \mathbf{x}\mathbf{x}^T\mathbf{M})\mathbf{x}_{\perp}' \\
    \label{eq:par-property}
    \mathbf{x}_{\parallel}' &= \mathbf{x}\mathbf{x}^T\mathbf{M}\mathbf{x}_{\parallel}'.
\end{align}
Substituting $\mathbf{x}'$ to its parallel and orthogonal components to equation \ref{eq:diff-single-norm} yields
\begin{equation}
    \mathbf{x}^T\mathbf{Mx}_{\parallel}' = -\frac{1}{2}\mathbf{x}^T\mathbf{M'x}.
\end{equation}
Note that the orthogonal component disappears from the equation above based on its property on equation \ref{eq:ortho-property}.
Multiplying the equation above with $\mathbf{x}$ from left on both sides, then use the property in equation \ref{eq:par-property} gives the parallel component of $\mathbf{x}'$,
\begin{equation}
    \label{eq:fwdderiv-xpar}
    \mathbf{x}_{\parallel}' = -\frac{1}{2} \mathbf{x}\mathbf{x}^T\mathbf{M'x}.
\end{equation}

The orthogonal component can be obtained by substituting $\mathbf{x}'$ to its component to equation \ref{eq:fwdderiv-before-splitting},
\begin{equation}
    (\mathbf{A} - \lambda\mathbf{M})\mathbf{x}_{\perp}' = -(\mathbf{I} - \mathbf{Mxx}^T)(\mathbf{A}' - \lambda\mathbf{M}')\mathbf{x}.
\end{equation}
Note that $(\mathbf{A} - \lambda\mathbf{M})\mathbf{x}_{\parallel}' = \mathbf{0}$ which can be shown using equation \ref{eq:par-property} and \ref{eq:single-eigdecomp}.
The equation above can be satisfied with finite $\mathbf{x}_\perp'$ because the vectors on both sides do not have a component parallel to $\mathbf{x}$.
Therefore, the orthogonal component of $\mathbf{x}'$ can be obtained by
\begin{equation}
    \label{eq:fwdderiv-xperp}
    \mathbf{x}_\perp' = -(\mathbf{I} - \mathbf{xx}^T\mathbf{M})(\mathbf{A} - \lambda\mathbf{M})^+(\mathbf{I} - \mathbf{Mxx}^T)(\mathbf{A}' - \lambda\mathbf{M}')\mathbf{x},
\end{equation}
where the plus superscript ($\cdot^+$) is the pseudo-inverse and
the term $(\mathbf{I} - \mathbf{xx}^T\mathbf{M})$ is applied to ensure the orthogonality of $\mathbf{x}_\perp'$ with respect to $\mathbf{x}$.
Combining the parallel and orthogonal components of $\mathbf{x}'$ gives the expression of forward derivative of the eigenvector,
\begin{equation}
    \label{eq:fwdderiv-single-eivec}
    \mathbf{x}' = -\frac{1}{2} \mathbf{x}\mathbf{x}^T\mathbf{M'x}-(\mathbf{I} - \mathbf{xx}^T\mathbf{M})(\mathbf{A} - \lambda\mathbf{M})^+(\mathbf{I} - \mathbf{Mxx}^T)(\mathbf{A}' - \lambda\mathbf{M}')\mathbf{x}.
\end{equation}

With the derivative of a single eigenvalue in equation \ref{eq:fwdderiv-single-eival} and eigenvector in equation \ref{eq:fwdderiv-single-eivec}, the expression for all $k$ eigenvalues and eigenvectors can be written as
\begin{align}
    \label{eq:fwdderiv-multi-eivals}
    \mathbf{\Lambda}' &= \mathbf{I}\circ \left[\mathbf{X}^T \left(\mathbf{A}' \mathbf{X} - \mathbf{M'X\Lambda}\right)\right] \\
    \mathbf{X}' &= -\frac{1}{2}\mathbf{X}\left[\mathbf{I}\circ\left(\mathbf{X}^T\mathbf{M'X}\right)\right] - \mathbf{Y}' + \mathbf{X}\left[\mathbf{I}\circ\left(\mathbf{X}^T\mathbf{MY'}\right)\right]
\end{align}
with $\circ$ denotes element-wise multiplication and
\begin{align}
    \label{eq:temp-y-fwdderiv}
    \mathbf{AY' - MY'\Lambda} &= \mathbf{V}' - \mathbf{MX}\left[\mathbf{I}\circ\left(\mathbf{X}^T\mathbf{V'}\right)\right] \\
    \label{eq:temp-v-fwdderiv}
    \mathbf{V'} &= \mathbf{A'X} - \mathbf{M'X\Lambda}.
\end{align}

\subsection{Backward derivative}

Once the forward derivative is found, the backward derivative can be found relatively easily.
If the forward derivative of a matrix can be expressed as
\begin{equation}
    \mathbf{P'} = \mathbf{QR'S},
\end{equation}
then by a simple manipulation using index notation, the backward derivative can be expressed as
\begin{equation}
    \label{eq:bckderiv-general-rule}
    \mathbf{\overline{R}} = \mathbf{Q}^T\mathbf{\overline{P}}\mathbf{S}^T.
\end{equation}

Using equations \ref{eq:bckderiv-general-rule}, \ref{eq:fwdderiv-single-eival}, and \ref{eq:fwdderiv-single-eivec}, the backward derivatives of $\mathbf{A}$ and $\mathbf{M}$ from one eigenvalue and eigenvector are
\begin{align}
    \label{eq:bckderiv-single-A}
    \mathbf{\overline{A}} &= \mathbf{xx}^T \overline{\lambda} -
(\mathbf{I} - \mathbf{xx}^T\mathbf{M})(\mathbf{A} - \lambda \mathbf{M})^{+}
(\mathbf{I} - \mathbf{Mxx}^T)\mathbf{\overline{x}} \mathbf{x}^T \\
\label{eq:bckderiv-single-M}
\mathbf{\overline{M}} &= -\mathbf{xx}^T \lambda \overline{\lambda}
-\frac{1}{2}\mathbf{xx}^T\mathbf{\overline{x}}\mathbf{x}^T +
\lambda (\mathbf{I} - \mathbf{xx}^T\mathbf{M})(\mathbf{A} - \lambda \mathbf{M})^{+}
(\mathbf{I} - \mathbf{Mxx}^T)\mathbf{\overline{x}} \mathbf{x}^T.
\end{align}
For the case with $k$ eigenvalues and eigenvectors, the contribution from every single pair must be summed.
Therefore, it can be written as
\begin{align}
    \mathbf{\overline{A}} &= \mathbf{X\overline{\Lambda}X}^T - \mathbf{\overline{V}X}^T\\
    \mathbf{\overline{M}} &= \mathbf{X\Lambda\overline{\Lambda}X}^T -
        \frac{1}{2}\mathbf{X}\left[\mathbf{I}\circ\left(\mathbf{X}^T\mathbf{\overline{X}}\right)\right]\mathbf{X}^T +
        \mathbf{\overline{V}\Lambda X}^T.
\end{align}
where $\circ$ denotes the element-wise multiplication and 
\begin{align}
    \label{eq:temp-y-bckderiv}
    \mathbf{\overline{V}} &= \mathbf{\overline{Y}} -
        \mathbf{X}\left[\mathbf{I}\circ\left(\mathbf{X}^T\mathbf{M\overline{Y}}\right)\right] \\
    \label{eq:temp-v-bckderiv}
    \mathbf{A\overline{Y}} - \mathbf{M\overline{Y}E} &=
        \mathbf{\overline{X}} -
        \mathbf{MX} \left[\mathbf{I}\circ\left(\mathbf{X}^T \mathbf{\overline{X}}\right)\right].
\end{align}

\section{Requirements for degenerate case}

Now consider a case where there are repeated eigenvalues, i.e. the degenerate case.
In the degenerate case, not all perturbations $\mathbf{A'}$ and $\mathbf{M'}$ can yield finite perturbations of $\mathbf{X'}$ and $\mathbf{\Lambda'}$.
Similar things also apply for the backward derivative where not all value's function of eigenvectors and eigenvalues can produce the finite backward derivative of $\mathbf{A}$ and $\mathbf{M}$.
In this section, the requirements to get finite forward and backward are derived.

\subsection{Forward derivative}

Let's denote the set of indices that have the same eigenvalue as the $j$-th eigenvalue as $\mathrm{d}(j)$, i.e.
\begin{equation}
    \mathrm{d}(j) = \left\{i\in\mathbb{Z}\ |\ i < n, \lambda_i = \lambda_j, i\neq j\right\}.
\end{equation}

Consider, only for this section, that all eigenvectors and eigenvalues of $\mathbf{A}$ are available.
The complete eigenvectors matrix is denoted as $\mathbf{U}\in\mathbb{R}^{n\times n}$ while the complete eigenvalues diagonal matrix is denoted as $\mathbf{E}\in\mathbb{R}^{n\times n}$.
Therefore, equation \ref{eq:orig-eidecomp} can be written as
\begin{equation}
    \label{eq:complete-eidecomp}
    \mathbf{AU} = \mathbf{MUE}.
\end{equation}
For a positive definite symmetric matrix $\mathbf{M}$, the equation above can be rewritten as a simple eigendecomposition for a real symmetric matrix,
\begin{align}
    \mathbf{M}^{-1/2}\mathbf{AM}^{-1/2}\mathbf{W} &= \mathbf{WE},
\end{align}
where $\mathbf{W} = \mathbf{M}^{1/2}\mathbf{U}$ and $\mathbf{M}^{1/2}\mathbf{M}^{1/2} = \mathbf{M}$.
As $\mathbf{M}^{-1/2}\mathbf{AM}^{-1/2}$ is a real symmetric matrix, $\mathbf{W}$ is unitary, therefore
\begin{align}
    \mathbf{U}^T\mathbf{MU} &= \mathbf{I} \\
    \label{eq:uut}
    \mathbf{UU}^T &= \mathbf{M}^{-1}.
\end{align}

Applying $\mathbf{U}^T$ from the right on both sides of equation \ref{eq:complete-eidecomp} and using equation \ref{eq:uut}, the matrix $(\mathbf{A} - \lambda\mathbf{M})$ and its pseudo-inverse for an eigenvalue $\lambda$ can be written as
\begin{align}
    (\mathbf{A} - \lambda \mathbf{M}) &= \mathbf{MU}(\mathbf{E} - \lambda\mathbf{I})\mathbf{U}^T\mathbf{M} \\
    \label{eq:inverse-a-min-lambda-m}
    (\mathbf{A} - \lambda \mathbf{M})^{+} &= \mathbf{M}^{-1}\mathbf{U}^{-T}(\mathbf{E} - \lambda \mathbf{I})^{+}\mathbf{U}^{-1}\mathbf{M}^{-1}.
\end{align}
Substituting $(\mathbf{A} - \lambda \mathbf{M})^{+}$ to the equation \ref{eq:fwdderiv-xperp} for the orthogonal perturbation component of an eigenvector $\mathbf{x}_j$ yields
\begin{equation}
    \mathbf{x}_{j\perp}' = \sum_{i\neq j} \frac{\mathbf{x}_i\mathbf{x}_i^T(\mathbf{A'} - \lambda_j\mathbf{M'})\mathbf{x}_j}{\lambda_i - \lambda_j}.
\end{equation}

From the equation above, it can be seen for all $i\in\mathrm{d}(j)$, $\lambda_i = \lambda_j$, and a division by zero occurs.
To make the forward derivative finite, the numerator must be zero whenever $\lambda_i = \lambda_j$.
Therefore, the requirements to get finite forward derivative,
\begin{equation}
    \mathbf{x}_i^T(\mathbf{A'} - \lambda_j\mathbf{M'})\mathbf{x}_j = 0\ \forall\ i\in\mathrm{d}(j).
\end{equation}
To express the equation above in terms of the retrieved eigenvectors and eigenvalues matrices, $\mathbf{X}$ and $\mathbf{\Lambda}$, it can be written as
\begin{equation}
    \label{eq:fwdderiv-degen-requirements}
    (\mathbf{D} - \mathbf{I})\circ\left[\mathbf{X}^T(\mathbf{A'X} - \mathbf{M'X\Lambda})\right] = \mathbf{0}.
\end{equation}

\subsection{Backward derivative}

Similar to the forward derivative, in backward derivative the problem of the degenerate case can be found in the contribution from the orthogonal component of the eigenvector.
The contribution from the orthogonal component of the eigenvector $\mathbf{x}_j$ to $\mathbf{\overline{A}}$ and $\mathbf{\overline{M}}$ on equations \ref{eq:bckderiv-single-A} and \ref{eq:bckderiv-single-M} are
\begin{align}
    \mathbf{\overline{A}}_{j\perp} &= -
    (\mathbf{I} - \mathbf{x}_j\mathbf{x}_j^T\mathbf{M})(\mathbf{A} - \lambda_j \mathbf{M})^{+}
    (\mathbf{I} - \mathbf{Mx}_j\mathbf{x}_j^T)\mathbf{\overline{x}}_j \mathbf{x}_j^T\\
    \mathbf{\overline{M}}_{j\perp} &= \lambda_j
    (\mathbf{I} - \mathbf{x}_j\mathbf{x}_j^T\mathbf{M})(\mathbf{A} - \lambda_j \mathbf{M})^{+}
    (\mathbf{I} - \mathbf{Mx}_j\mathbf{x}_j^T)\mathbf{\overline{x}}_j \mathbf{x}_j^T.
\end{align}
Substituting $(\mathbf{A} - \lambda_j\mathbf{M})^{+}$ from equation \ref{eq:inverse-a-min-lambda-m} to the equations above produces
\begin{align}
    \mathbf{\overline{A}}_{j\perp} &= -\sum_{i\neq j}\frac{\mathbf{x}_i\mathbf{x}_i^T\mathbf{\overline{x}}_j\mathbf{x}_j^T}{\lambda_i - \lambda_j}\\
    \mathbf{\overline{M}}_{j\perp} &= \lambda_j \sum_{i\neq j}\frac{\mathbf{x}_i\mathbf{x}_i^T\mathbf{\overline{x}}_j\mathbf{x}_j^T}{\lambda_i - \lambda_j}.
\end{align}
The total contribution from all orthogonal components from all indices is just simply sum of the terms above for all $j$,
\begin{align}
    \label{eq:bckderiv-Aperp-contrib-series}
    \mathbf{\overline{A}}_{\perp} &= -\sum_j\sum_{i\neq j}\frac{\mathbf{x}_i\mathbf{x}_i^T\mathbf{\overline{x}}_j\mathbf{x}_j^T}{\lambda_i - \lambda_j}\\
    \label{eq:bckderiv-Mperp-contrib-series}
    \mathbf{\overline{M}}_{\perp} &= \sum_j \sum_{i\neq j}\lambda_j\frac{\mathbf{x}_i\mathbf{x}_i^T\mathbf{\overline{x}}_j\mathbf{x}_j^T}{\lambda_i - \lambda_j}.
\end{align}

If $\mathbf{A}$ and $\mathbf{M}$ are parameterized by $\theta_A$ and $\theta_M$ respectively, the backward derivative of $\theta_A$ and $\theta_M$ from the orthogonal eigenvector components are
\begin{align}
    \overline{\theta_A}_\perp &= \mathrm{tr}\left[\mathbf{\overline{A}}^T_\perp \frac{\partial \mathbf{A}}{\partial \theta_A}\right] \\
    \overline{\theta_M}_\perp &= \mathrm{tr}\left[\mathbf{\overline{M}}^T_\perp \frac{\partial \mathbf{M}}{\partial \theta_M}\right].
\end{align}
Substituting the equations above to equations \ref{eq:bckderiv-Aperp-contrib-series} and \ref{eq:bckderiv-Mperp-contrib-series}, then using the cyclic property of trace gives the expressions below,
\begin{align}
    \overline{\theta_A}_\perp &= -\sum_j\sum_{i\neq j} \frac{1}{\lambda_i - \lambda_j}\left(\mathbf{x}_i^T\mathbf{\overline{x}}_j\mathbf{x}_i^T\frac{\partial \mathbf{A}}{\partial \theta_A}\mathbf{x}_j\right) \\
    \overline{\theta_M}_\perp &= \sum_j\sum_{i\neq j} \frac{\lambda_j}{\lambda_i - \lambda_j}\left(\mathbf{x}_i^T\mathbf{\overline{x}}_j\mathbf{x}_i^T\frac{\partial \mathbf{M}}{\partial \theta_M}\mathbf{x}_j\right).
\end{align}

Again, if $\lambda_i = \lambda_j$ for $i \in \mathrm{d}(j)$, the denumerator becomes 0.
To solve this problem, the term where $\lambda_i = \lambda_j$ must be vanished.
For two indices, $i$ and $j$ where $i \in \mathrm{d}(j)$, the terms involving $i$ and $j$ in the equations above are
\begin{align}
    \overline{\theta_A}_\perp &= -\frac{1}{\lambda_i - \lambda_j}\left(\mathbf{x}_i^T\mathbf{\overline{x}}_j\mathbf{x}_i^T\frac{\partial \mathbf{A}}{\partial \theta_A}\mathbf{x}_j - \mathbf{x}_j^T\mathbf{\overline{x}}_i\mathbf{x}_j^T\frac{\partial \mathbf{A}}{\partial \theta_A}\mathbf{x}_i\right) + ... \\
    \overline{\theta_M}_\perp &= \frac{\lambda_j}{\lambda_i - \lambda_j}\left(\mathbf{x}_i^T\mathbf{\overline{x}}_j\mathbf{x}_i^T\frac{\partial \mathbf{M}}{\partial \theta_M}\mathbf{x}_j - \mathbf{x}_j^T\mathbf{\overline{x}}_i\mathbf{x}_j^T\frac{\partial \mathbf{M}}{\partial \theta_M}\mathbf{x}_i\right) + ...
\end{align}
To make the terms vanish, the conditions below must be satisfied,
\begin{align}
    \mathbf{x}_i^T\mathbf{\overline{x}}_j\mathbf{x}_i^T\frac{\partial \mathbf{A}}{\partial \theta_A}\mathbf{x}_j &= \mathbf{x}_j^T\mathbf{\overline{x}}_i\mathbf{x}_j^T\frac{\partial \mathbf{A}}{\partial \theta_A}\mathbf{x}_i \\
    \lambda_j\mathbf{x}_i^T\mathbf{\overline{x}}_j\mathbf{x}_i^T\frac{\partial \mathbf{M}}{\partial \theta_M}\mathbf{x}_j &= \lambda_j\mathbf{x}_j^T\mathbf{\overline{x}}_i\mathbf{x}_j^T\frac{\partial \mathbf{M}}{\partial \theta_M}\mathbf{x}_i.
\end{align}
If $\mathbf{A}$ and $\mathbf{M}$ are always symmetric, i.e. $\partial \mathbf{A}/\partial \theta_A$ and $\partial \mathbf{M} / \partial \theta_M$ are symmetric, the conditions to get finite backward derivative for $\mathbf{A}$ and $\mathbf{M}$ become
\begin{equation}
    \mathbf{x}_i^T\mathbf{\overline{x}}_j = \mathbf{x}_j^T\mathbf{\overline{x}}_i\ \forall\ i\in\mathrm{d}(j)
\end{equation}
or
\begin{align}
    \mathbf{x}_i^T\frac{\partial\mathbf{A}}{\partial \theta_A}\mathbf{x}_j &= 0\ \forall\ i\in\mathrm{d}(j) \\
    \lambda_j\mathbf{x}_i^T\frac{\partial\mathbf{M}}{\partial \theta_M}\mathbf{x}_j &= 0\ \forall\ i\in\mathrm{d}(j).
\end{align}

In terms of multiple eigenvalues and eigenvectors matrices, the requirements above can be written as
\begin{align}
    (\mathbf{D} - \mathbf{I})\circ\left(\mathbf{X}^T\mathbf{A'X}\right) &= \mathbf{0} \\
    (\mathbf{D} - \mathbf{I})\circ\left(\mathbf{X}^T\mathbf{M'X\Lambda}\right) &= \mathbf{0}
\end{align}
and
\begin{equation}
    (\mathbf{D} - \mathbf{I}) \circ \left(\mathbf{X}^T\mathbf{\overline{X}} - \mathbf{\overline{X}}^T\mathbf{X}\right) = \mathbf{0}
\end{equation}

\section{Derivatives for the degenerate case}
\subsection{Forward derivative}
If the conditions in the previous section are satisfied, then the forward and backward derivatives can be calculated.
Let's denote the set of eigenvectors that have the same corresponding eigenvalue equals to $\lambda$ as $\mathrm{g}(\lambda)$, i.e.
\begin{equation}
    \mathrm{g}(\lambda) = \left\{\mathbf{x}_i\in\mathbb{R}^{n\times 1}\ |\ \mathbf{Ax}_i = \lambda_i\mathbf{Mx}_i,\ \lambda_i = \lambda\right\},
\end{equation}
and the matrix containing the eigenvectors $\mathrm{g}(\lambda)$ as its columns as $\mathbf{X}_{\mathrm{g}(\lambda)}$.

In the degenerate case the forward derivative of the eigenvalues in equations \ref{eq:fwdderiv-single-eival} and \ref{eq:fwdderiv-multi-eivals} are still valid.
However, the forward derivative of the eigenvector from the previous section is invalid because the matrix $(\mathbf{A} - \lambda\mathbf{M})$ in equation \ref{eq:fwdderiv-before-splitting} have the rank of $n - |\mathrm{g}(\lambda)|$ instead of $n-1$.

To derive the perturbation of the eigenvector in the degenerate case, let's rewrite the equation \ref{eq:fwdderiv-before-splitting}
\begin{equation}
    \label{eq:diff-single-eigdecomp2}
    (\mathbf{A} - \lambda\mathbf{M})\mathbf{x}' = -(\mathbf{I} - \mathbf{Mxx}^T)(\mathbf{A}' - \lambda\mathbf{M}')\mathbf{x}.
\end{equation}
The matrix $(\mathbf{A} - \lambda\mathbf{M})$ nullify any vector components that are parallel to any eigenvectors in $\mathrm{g}(\lambda)$.
Therefore, $\mathbf{x}'$ needs to be treated separately based on its components : (1) the parallel component, $\mathbf{x}_\parallel'$, that is parallel to $\mathbf{x}$, (2) the degenerate component, $\mathbf{x}_\circ'$, that can be expressed by a linear combination of eigenvectors in $\left[\mathrm{g}(\lambda) - \{\mathbf{x}\}\right]$, and (3) the perpendicular component, $\mathbf{x}_\perp'$, that is perpendicular in $\mathbf{M}$ to the eigenvectors in $\mathrm{g}(\lambda)$,
\begin{equation}
    \mathbf{x}' = \mathbf{x}_\parallel' + \mathbf{x}_\circ' + \mathbf{x}_\perp'.
\end{equation}
The properties of those components are
\begin{align}
    \mathbf{x}_\parallel' &= \mathbf{xx}^T\mathbf{Mx}_\parallel' \\
    \mathbf{x}_\circ' &= \left(\mathbf{X}_{\mathrm{g}(\lambda)}\mathbf{X}_{\mathrm{g}(\lambda)}^T - \mathbf{xx}^T\right)\mathbf{M}\mathbf{x}_\circ'\\
    \mathbf{x}_\perp' &= \left(\mathbf{I} - \mathbf{X}_{\mathrm{g}(\lambda)}\mathbf{X}_{\mathrm{g}(\lambda)}^T\mathbf{M}\right)\mathbf{x}_\perp'.
\end{align}

Splitting $\mathbf{x}'$ into its components to equation \ref{eq:diff-single-norm} will eliminate $\mathbf{x}_\circ'$ and $\mathbf{x}_\perp'$, and only $\mathbf{x}_\parallel'$ remains.
By following the same steps as in the non-degenerate case, the parallel component of the eigenvector's perturbation in the degenerate case has the same expression as in the non-degenerate case,
\begin{equation}
    \mathbf{x}_\parallel' = -\frac{1}{2}\mathbf{xx}^T\mathbf{M}'\mathbf{x}.
\end{equation}

Now the task is to find the orthogonal component.
Substituting $\mathbf{x}'$ to its 3 components into equation \ref{eq:diff-single-eigdecomp} will eliminate the degenerate and parallel components, $\mathbf{x}_\circ'$ and $\mathbf{x}_\parallel'$.
Therefore, it can be written as
\begin{equation}
    \label{eq:degen-bckderiv-before-pseudoinv}
    (\mathbf{A} - \lambda\mathbf{M})\mathbf{x}_\perp' = -(\mathbf{I} - \mathbf{Mxx}^T)(\mathbf{A}' - \lambda\mathbf{M}')\mathbf{x}.
\end{equation}
Notice that the vector $(\mathbf{A} - \lambda\mathbf{M})\mathbf{x}_\perp'$ on the left hand side has no components parallel to the eigenvectors in $\mathrm{g}(\lambda)$, i.e.
\begin{equation}
    \mathbf{X}_{\mathrm{g}(\lambda)}^T(\mathbf{A} - \lambda\mathbf{M})\mathbf{x}_\perp' = \mathbf{0}.
\end{equation}
Therefore, to get finite $\mathbf{x}_\perp'$, the vector on the right hand side must have 0 components parallel to the columns of $\mathbf{X}_{\mathrm{g}(\lambda)}$, i.e.
\begin{equation}
    \mathbf{X}_{\mathrm{g}(\lambda)}^T(\mathbf{I} - \mathbf{Mxx}^T)(\mathbf{A}' - \lambda\mathbf{M}')\mathbf{x} = \mathbf{0}.
\end{equation}

The condition above is satisfied if \ref{eq:fwdderiv-degen-requirements} satisfied.
Applying the pseudo-inverse of $\left(\mathbf{A} - \lambda\mathbf{M}\right)$ to the equation \ref{eq:degen-bckderiv-before-pseudoinv} and applying the orthogonalization operator gives
\begin{equation}
    \mathbf{x}_\perp' = -\left(\mathbf{I} - \mathbf{X}_{\mathrm{g}(\lambda)}\mathbf{X}_{\mathrm{g}(\lambda)}^T\mathbf{M}\right)\left(\mathbf{A} - \lambda\mathbf{M}\right)^+\left(\mathbf{I} - \mathbf{Mxx}^T\right)\left(\mathbf{A'} - \lambda\mathbf{M'}\right)\mathbf{x}
\end{equation}
where the term $\left(\mathbf{I} - \mathbf{X}_{\mathrm{g}(\lambda)}\mathbf{X}_{\mathrm{g}(\lambda)}^T\mathbf{M}\right)$ is applied to ensure the orthogonality of $\mathbf{x}_\perp'$ with respect to all columns of $\mathbf{X}_{\mathrm{g}(\lambda)}$.
The expression is similar to equation \ref{eq:fwdderiv-xperp} of the perpendicular component in the non-degenerate case, except for the orthogonalization term.

Having found the parallel and perpendicular components, the degenerate component, $\mathbf{x}_\circ'$ still remains to be found.
As the degenerate component is eliminated in equation \ref{eq:diff-single-norm} and \ref{eq:diff-single-eigdecomp2}, there is no restriction applies to the degenerate component and thus it can take any value.
To make it simple, the degenerate component can be assigned to be zeros,
\begin{equation}
    \mathbf{x}_\circ' = \mathbf{0}.
\end{equation}
Combining the results from its components, the eigenvector perturbation can be written as,
\begin{equation}
    \label{eq:fwdderiv-single-eivec-degen}
    \mathbf{x}' = -\frac{1}{2}\mathbf{xx}^T\mathbf{M}'\mathbf{x} -\left(\mathbf{I} - \mathbf{X}_{\mathrm{g}(\lambda)}\mathbf{X}_{\mathrm{g}(\lambda)}^T\mathbf{M}\right)\left(\mathbf{A} - \lambda\mathbf{M}\right)^+\left(\mathbf{I} - \mathbf{Mxx}^T\right)\left(\mathbf{A'} - \lambda\mathbf{M'}\right)\mathbf{x}.
\end{equation}

With the perturbation of a single eigenvalue and eigenvector found in the equations above, we can write the expression for the perturbation of $k$ eigenvalues and eigenvectors,
\begin{align}
    \mathbf{\Lambda}' &= \mathbf{I}\circ \left[\mathbf{X}^T \left(\mathbf{A}' \mathbf{X} - \mathbf{M'X\Lambda}\right)\right] \\
    \mathbf{X}' &= -\frac{1}{2}\mathbf{X}\left[\mathbf{I}\circ\left(\mathbf{X}^T\mathbf{M'X}\right)\right] - \mathbf{Y}' + \mathbf{X}\left[\mathbf{D}\circ\left(\mathbf{X}^T\mathbf{MY'}\right)\right]
\end{align}
where $\mathbf{Y'}$ and $\mathbf{V'}$ are given in equations \ref{eq:temp-y-fwdderiv} and \ref{eq:temp-v-fwdderiv}, respectively, and $\mathbf{D}$ is the degenerate matrix with elements
\begin{equation}
    \label{eq:degen-matrix}
    D_{ij} = \begin{cases}
        1, &\ \text{if }\Lambda_{ii} = \Lambda_{jj}\\
        0, &\ \text{otherwise}.
    \end{cases}
\end{equation}

\subsection{Backward derivative}

Following the same process in obtaining the backward derivative from the forward derivatives in equations \ref{eq:fwdderiv-single-eivec-degen} and \ref{eq:fwdderiv-single-eival} produces
\begin{align}
    \mathbf{\overline{A}} &= \mathbf{xx}^T \overline{\lambda} -
\left(\mathbf{I} - \mathbf{xx}^T\mathbf{M}\right)(\mathbf{A} - \lambda \mathbf{M})^{+}
\left(\mathbf{I} - \mathbf{MX}_{\mathrm{g}(\lambda)}\mathbf{X}_{\mathrm{g}(\lambda)}^T\right)\mathbf{\overline{x}} \mathbf{x}^T \\
\mathbf{\overline{M}} &= -\mathbf{xx}^T \lambda \overline{\lambda}
-\frac{1}{2}\mathbf{xx}^T\mathbf{\overline{x}}\mathbf{x}^T +
\lambda \left(\mathbf{I} - \mathbf{xx}^T\mathbf{M}\right)(\mathbf{A} - \lambda \mathbf{M})^{+}
\left(\mathbf{I} - \mathbf{MX}_{\mathrm{g}(\lambda)}\mathbf{X}_{\mathrm{g}(\lambda)}^T\right)\mathbf{\overline{x}} \mathbf{x}^T.
\end{align}
The expressions are similar to the non-degenerate case, except for the term $\left(\mathbf{I} - \mathbf{MX}_{\mathrm{g}(\lambda)}\mathbf{X}_{\mathrm{g}(\lambda)}^T\right)$ on the right hand side of the equation.

For multiple eigenvalues and eigenvectors, the expression reads
\begin{align}
    \mathbf{\overline{A}} &= \mathbf{X\overline{\Lambda}X}^T - \mathbf{\overline{V}X}^T\\
    \mathbf{\overline{M}} &= \mathbf{X\Lambda\overline{\Lambda}X}^T -
        \frac{1}{2}\mathbf{X}\left[\mathbf{I}\circ\left(\mathbf{X}^T\mathbf{\overline{X}}\right)\right]\mathbf{X}^T +
        \mathbf{\overline{V}\Lambda X}^T.
\end{align}
where
\begin{align}
    \label{eq:temp-y-bckderiv-degen}
    \mathbf{\overline{V}} &= \mathbf{\overline{Y}} -
        \mathbf{X}\left[\mathbf{I}\circ\left(\mathbf{X}^T\mathbf{M\overline{Y}}\right)\right] \\
    \label{eq:temp-v-bckderiv-degen}
    \mathbf{A\overline{Y}} - \mathbf{M\overline{Y}E} &=
        \mathbf{\overline{X}} -
        \mathbf{MX} \left[\mathbf{D}\circ\left(\mathbf{X}^T \mathbf{\overline{X}}\right)\right].
\end{align}
The expressions are very similar to the non-degenerate case, except on the equation \ref{eq:temp-v-bckderiv-degen} where the element-wise multiplication with an identity matrix is now replaced by the degenerate matrix $\mathbf{D}$.

\bigskip
{\bf Acknowledgment.} M.F.K. acknowledge support from the U.K. EPSRC under Grant No. EP/P015794/1.

\printbibliography

\end{document}